\documentclass[article,1p]{elsarticle}
\usepackage[english]{babel}
\usepackage{caption,longtable}
\usepackage{algorithm,algorithmic}
\usepackage{graphicx,color,float,rotating}
\usepackage{amssymb,amsmath}
\usepackage{epstopdf}
\newtheorem{theorem}{Theorem}

\newtheorem{corollary}{Corollary}
\newtheorem{definition}{Definition}
\newcommand{\be}{\begin{equation}}
\newcommand{\ee}{\end{equation}}

\begin{document}

\begin{frontmatter}

\title{On the Adimensional Scale Invariant Steffensen (ASIS) Method.}%

\author{ Vicente F. Candela\fnref{vc}}
\ead{vicente.candela@uv.es}
\address{Departament de Matem\`atiques, Universitat de Val\`encia.\\
e-mail: vicente.candela@uv.es}
\fntext[vc]{ESI International Chair@CEU-UCH. Departament de
Matem`atiques, Campus de Burjassot, Universidad de Valencia Carrer del Dr.
Moliner 50, 46100 Burjassot (Valencia), Spain.}

\begin{abstract}
Dimensionality of parameters and variables  is a fundamental issue in physics but mostly ignored from a mathematical point of view. Difficulties arising from dimensional inconsistence are overcome by scaling analysis and, often, both concepts, dimensionality and scaling, are confused. In the particular case of iterative methods for solving nonlinear equations, dimensionality and scaling affects their robutness: while some classical methods, such as Newton, are adimensional and scale independent, some other iterations as Steffensen's are not; their convergence depends on the scaling, and their evaluation needs a dimensional congruence. In this paper we introduce the concept of {\it adimensional form} of a function in order to study the behavior of iterative methods, thus correcting, if possible, some pathological features. From this adimensional form we will devise an adimensional and scale invariant method based on Steffensen's which we will call ASIS method. 
\end{abstract}


\end{frontmatter}

\section{Introduction}
The main difference between absolute and relative errors in approximation is that, while absolute errors share the dimension of the values they approximate, relative errors are adimensional, and therefore, independent from scaling. This trivial idea is often ignored when the parameter we want to approximate is the solution of an equation. From a numerical point of view, scaling is a fundamental feature, but,  often, dimensional congruence is ignored. As a consequence of disregarding dimensionality, some methods  devised originally for scalar, real or complex, equations
are not useful, or hard to adapt, to general vector, Banach spaces. This is a reason why many iterative methods are not used in optimization (of course, storage and computational cost are also important).

In \cite{meinsma}, a nice introduction to dimensional and scaling analysis is found. Although dimensional analysis goes back as far as Fourier in 1822, it was boosted in the second decade of last century, especially since the works of E. Buckingham (\cite{buck}, for example), most of them from a physical rather than mathematical point of view. In fact, dimensionality is a fundamental concept in physics in order to determine the consistence of equations and magnitudes, but not so important in mathematics. The simple example of coefficients of polynomials illustrates this assertion.

Given $p(x)=a_0+a_1x+ a_2x^2+ \cdots +a_nx^n$, if $x$ and $p(x)$ have dimensions denoted by $[x]$ and $[p]$ respectively, then the coefficients must have dimension $[a_k]=[p][x]^{-k}=[p^{(k)}(x)]$. This obvious remark implies that $a_k$ is a constant magnitude, but not a constant number (a change of scale both in $x$ or in $p(x)$ modifies its value). Numerically, this means that we must be careful with scaling (if $x$ is a length, a change from centimeters to inches, for instance, implies a change of coefficients), and one cannot forget the fact that $a_k$ is the $k$-th derivative of $p(x)$. Dimensional considerations do not come from the physical units we use, but from the kind of magnitude we are measuring (let's say time, length, weight, ...).  To put it shortly, dimensionality is inherent to the variable while scaling is accidental and subject to changes.

\vskip 6pt

In general, dimensional analysis spans over different, and apparently unrelated, topics in mathematics. It reaches its strength in nonlinear equations, where dimensionality and scaling are clearly separated. The purpose of this paper is the application of dimensional analysis of iterative methods for nonlinear equations. 

We start by some basic notions about convergent sequences. It is well known \cite{potra} that a convergent sequence, $\{x_k\}_{k \geq 0}$, in a Banach space $(E, || \cdot ||)$ with limit $x^* \in E$ has Q-order of convergence at least $p$ if there exists a constant $K>0$ such that:
\begin{equation}\label{qorder}
\frac{||x_{n+1}-x^*||}{||x_n-x^*||^p}\leq K
\end{equation}
\hskip 12pt Therefore, $K$ is not an adimensional parameter when $p \neq 1$. It has dimension $[K]=[x]^{1-p}$. Any change of scale $y=cx$ changes the value of $K$. This fact, that could become dangerous is no  more than annoying because there  exists another weaker definition, R-order of convergence \cite{potra}, which is scale invariant, where the sequence has R-order at least $p>1$ if there exists a constant $K$ and an $n_0>0$ such that:
\begin{equation}\label{rorder}
\limsup_{n\longrightarrow \infty} ||x_n-x^*||^{1/p^n} \leq K
\end{equation}
\hskip 12pt In this definition, $K$ is adimensional ($[K]=[1]$), and there are not issues derived from scales or dimensions. Besides the evident differences between both definitions, there is a conceptual motivation: Q-order controls the evolution of every term $n$, while R-order focus on the global (limit) behavior of the sequence. From this point of view, it is not surprising that Q-order is a more restrictive definition than that of R-order. 

Here, we need to establish an order independent from dimensionality, as it happens with R-order, but able to analyze convergence for any step $n$ as Q-order. Thus, we will introduce the {\it adimensional Q-order}, scale independent and adimensional, weaker than Q-order but keeping its essential core. As a consequence, we will get the {\it adimensional} R-order, which is, in fact, a reformulation of the traditional R-order.

\vskip 6pt

Now we are ready to focus on iterative methods for solving nonlinear equations as $F(x)=0$, when $F: E \longrightarrow V$, $E$ and $V$ Banach spaces. The idea behind these methods is correcting a previous approximation $x_k$  to the solution $x^*$ by:
\begin{equation}\label{delta}
x_{k+1}=x_k+\Delta_k=x_k \cdot (1+\delta_k)
\end{equation} 
$[\Delta_k]=[x]$. Therefore, $\Delta_k$, which is the absolute error, depends on $x$: any change of scale or dimension in $x$ leads to a similar change on $\Delta_k$. 
The relative error, $\delta_k$, is adimensional, and $\delta_k=\Delta_k /x_k$, relates both definitions. As we said above, it is not surprising that, numerically, the first option is preferred, though $\delta_k$ is free from any changes of magnitude.

\vskip 6pt

One of the simplest and best known methods is the first order one  obtained when $\Delta_k=-c F(x_k)$ under certain restrictions on $c$. Numerically, these restrictions are due to the scale: $||I-cF'(x)||<1$ for  $x$ in the domain. From a dimensional point of view, 
$$
[\Delta_k]=[c][F] \ \implies [c]=[x][F]^{-1}
$$
Restrictions on the size of $c$, are a consequence of its dimension. Considering $c$ dependent of $F'$ is more clarifying than thinking of $c$ as a constant.

On the other hand,   $\Delta_k=-F'(x_k)^{-1}F(x_k)$  gives place to Newton's, which is dimensionally consistent:

$
[\Delta_k]=[F][F'^{-1}]=[x]
$
\noindent (because $[F']=[F][x]^{-1}$).

\vskip 6pt

Another second order method such as Steffensen's, is defined as:
$$
\Delta_k=-F[x_k+F(x_k),x_k]^{-1}F(x_k)
$$
\noindent where $F[\cdot , \cdot ]$ is the divided difference operator (bracket notation should not be confused with the dimension). Scaling is indeed a concern, but the most important feature is the inconsistence of dimensions. The sum $x+F(x)$ in the operator  makes sense only if $E=V$, because $[x]$ and $[F]$ are not equal in general. This is one of the reasons why, though competitive with respect to Newton's, Steffensen's is less popular.

However, Steffensen's has a clear advantage over Newton's method: it is not necessary to evaluate Jacobians or gradients for Steffensen's, because divided differences are usually obtained by the interpolatory condition:
$$
F[x,y](x-y)=F(x)-F(y)
$$
\noindent which does not require derivatives. 

\vskip 6pt

The rest of this study will  focus on two main goals. The first one is the concept of {\it adimensional form} of a function, a tool we introduce which keeps all the information of the function, but it is invariant to changes of scales and dimensions of the variables. We will use these adimensional functions in a theoretical way, to obtain semilocal convergent results, and, in practical applications, to get versions of non robust methods in order to make them scale independent.  The second part of this paper deals with the analysis, development and features of these modified methods.

In particular, we will prove the power of adimensional functions in the particular case of Steffensen's method. We will devise an adimensional and scale invariant (ASIS) correction of Steffensen's method and we will settle some sufficient semilocal conditions of convergence, from which optimal estimates of error will be obtained. As a difference with previous published results, sufficient conditions of convergence will not be more restrictive than those of Newton's. Furthermore, so far, optimality has not been obtained for the Steffensen's method in literature; thus our estimates improve previous works in some aspects. 

\vskip 6pt

This paper is organized as follows: the next section will be divided into two subsections in which we review some basic and classic results on {\it dimensionality} and iterative methods for nonlinear equations in Banach spaces; we also point out there a fundamental theorem relating both worlds (dimensionality and  iterative methods). In \S\ref{sec3} we define the {\it adimensional form} of polynomials and functions and we will apply them to the analysis of semilocal convergence of iterative methods in \S\ref{sec4}. In \S\ref{sec5}  we devise new variations of classic scale dependent iterative methods, such as Steffensen's, in order to eliminate their scale and dimensional dependence. \S\ref{sec6} will illustrate through examples the theoretical results in this paper. Finally, we will draw conclusions in \S\ref{sec7}.

\section{Basic concepts on dimensionality and iterative methods}\label{sec2}

Dimension of a variable is a physical concept related to the units in which that variable is measured. It is one of those concepts that are easier to understand than to define. A physical quantity may represent time, length, weight or any other magnitude. Equations in physics relate one or more of these dimensions and they must show consistence: dimensions in both sides of the  equation must be equal. Nevertheless, sometimes this trivial remark is ignored because either it is not so evident to get the dimension of any term or, as it happens in mathematics, good scaling may overcome the potential risks of dimensionality.

But dimensionality and scale are not always equivalent. For example, changing the measure of a length from, say, centimeters to inches is a change of scale (and very dangerous indeed), but their dimensions are the same: length. Inconsistence happens when there are different measures involved in each term of the equation (time versus length, for instance). The problem increases when the equation is not scalar but defined in Banach spaces. In this case it is easy to confuse dimensionality and dimension (two vectors of the same dimension may belong to different Banach spaces and have different dimensionality even if there exists an isomorphism between both), leading to some problems in areas such as optimization.

In the rest of this section we introduce some aspects of dimensional analysis in \S\ref{2.1}, and a basic theory on iterative methods for nonlinear equations in \S\ref{2.2}.

\subsection{Dimensionality}\label{2.1}

From now on, we will use the classical notation for the dimension of a variable $x$, $[x]$. If  $z$ is an adimensional parameter, we will denote it by $[z]=[1]$. 
Unlike the classic dot notation in physics, we will denote the derivatives of a function $f(x)$ by $f'(x)$, $f''(x)$, ...  as usual in mathematics. We will omit the variables unless we shall explicit them due to the context.

For vectors and Banach spaces in general, the functions will be denoted by capital letters ($F(x)$, $G(x)$, ...), and the vectors $x$, $y$, ... will use the same notation as for scalars (and there shall not be any confusion). $F'(x)$ and $F''(x)$ will denote the Jacobian and Hessian of $F(x)$, respectively.

$F: E \longrightarrow V$, where $E$ and $V$ are Banach spaces. Their norms $(E, || \cdot ||)$, and $(V, || \cdot ||)$, though different, will be denoted with the same notation (and once again hoping it shall not be confusing).

It is not intuitive to introduce dimensionality in nonscalar Banach spaces. When the space is finite dimensional, we define $[x]$ as a vector where every component is the dimension of the corresponding element of the base. Thus, if $x=(x_1, \cdots , x_m)$, then $[x]=([x_1], \cdots , [x_m])$. We must consider other operations besides the sum and multiplication by scalars in the context of nonlinear equations. Hence, nonlinear operations must be carried pointwise (component by component). For example, $x \ {\rm o} \  y=(x_1 {\rm o}  \ y_1 , \cdots , 
x_m {\rm o} \ y_m)$ for any operation ${\rm o}$.

For infinite dimension Banach spaces, dimensionality loses its physical sense and it is a more complicated concept. For the purpose of this paper, we consider dimensionality in a similar way as finite dimensional vector spaces: here, $[x]$ is an infinite vector, and its operations are also done pointwise. There is not misunderstanding when all the bases have the same dimensionality but, in general, the dimensional vectors may not be homogeneous.
  
Elementary rules of dimensionality apply:
$$
[x \cdot y]=[x][y] \ ; \ [x/y]=[x][y]^{-1} \ ; \ [x^n]=[x]^n
$$

Rules for differentiability and integration are particularly interesting:
$$
[F']=[F]/[x]=[F][x]^{-1} \ ; \
\left[\int F(x) {\rm d}x\right]=[F][x]
$$
\hskip 12pt Before going on, we will review some basics on iterative methods. This classic topic has been treated in most basic books on iterative methods (\cite{ost,ort}, for instance).

\subsection{Iterative methods for nonlinear equations}\label{2.2}

Iterative methods for solving a nonlinear scalar equation $f(x)=0$ are based on a simple idea: we start by approximating the value $x^*$ such that $f(x^*)=0$ by an initial guess $x_0$, and, after that, correcting the approximations by recurrence: $x_{k+1}=x_k+\Delta_k$. Thus, a sequence $\{x_n\}_{n\geq 0}$ is obtained. If that sequence converges to the solution $x^*$, we get the solution of the equation. 

Of course, this idea becomes more and more complicated in practice. In order to get convergence, the corrections $\Delta_k$ must be well devised (depending on the function $f(x)$ and, often, its derivatives if they exist), and some restrictions must apply on the own function and the initial guess. On the other side, a stopping criterion is needed in order to get a good approximation to the root.

The accuracy of the method depends not only on the design (consistence) of the recurrence but on the stability (sensitivity to errors) and the speed of convergence. There are different ways to define the speed of a sequence, but we point out Q-order and R-order defined in (\ref{qorder}) and (\ref{rorder}) respectively.

Perhaps the easiest iterative method is the {\it bisection}. Its approximates are not values but intervals. However, a point version of this method can be obtained if we consider the middle point of the interval where the function $f(x)$ switches signs. As the intervals are divided every step, a choice (an {\it if}) must be done every step. No other information about $f(x)$ is needed but the sign. When there is a change of sign, it has first order convergence:
$$
||x_{n+1}-x^*|| \approx \frac{1}{2} ||x_n-x^*||
$$
\hskip 12pt The more we know something about $f(x)$, the more sophisticated methods can be obtained. By carefully choosing a constant $c \in {\mathbb R}$, the recurrence  
\begin{equation}\label{firstorder}
\Delta_k=-c \cdot f(x_k)
\end{equation}
\noindent is a first order method.

As we saw in the introduction, we can eliminate dimensionality of the parameter $c$ by introducing this first order method:

\begin{equation}\label{alter}
x_{n+1}=x_n-\frac{\lambda}{f'(x_0)}f(x_k)
\end{equation}
$\lambda$ is adimensional, but some scale requirements are still needed: $0<|\lambda f'(x)/f'(x_0)|<2$. Improvement comes from the fact that, once $\lambda$ is fixed, any rescaling on $x$ or $f(x)$ does not affect $\lambda$.

However, the main method around which all the theory has evolved is Newton's, where:

\begin{equation}\label{newton}
\Delta_k= -\frac{f(x_k)}{f'(x_k)}
\end{equation}
\hskip 12pt It needs more requirements than the other two ($f(x)$ must be differentiable and, its derivative, invertible), and its convergence is usually local, but, when it works, it is a second order method. More important, it provides insight on this topic. A great amount of papers have been published on Newton's. Since Kantorovich settled conditions for convergence of this method in Banach spaces (\cite{kant}), lots of papers have appeared analyzing this method, some of them classic (\cite{den,tapia,gragg} among others) while new and different kind of variants and analysis are also frequent, \cite{ozban,weerakon,ezq,arg},  to cite just a few. In \cite{ezq} an updated bibliography can also be found).

In the cited bibliography some variants of Newton provided in literature can be found: quasi-Newton in optimization \cite{nocedal}, approximated Newton in computational analysis are some of these versions \cite{dembo,ye}.

Let us focus on quasi-Newton methods.  They are based on the approximation of the derivative by means of divided differences. For scalar functions, a divided difference operator in two nodes $x$, $y$, denoted by $f[x,y]$ (and, once again, we hope its context makes it clear when the brackets refer to dimension or to this operator) is defined as
\begin{equation}\label{divdif}
f[x,y]=\frac{f(x)-f(y)}{x-y}
\end{equation}
\hskip 12pt Considered as an operator, $f[x,y]$ verifies the interpolatory condition:
$$
f[x,y](x-y)=f(x)-f(y)
$$
\hskip 12pt If $f(x)$ is at least twice differentiable, divided differences approximate the derivative of $f(x)$:
$$
f[x,y]=f'(x)+f''(x+\xi (y-x))(y-x)
$$
\noindent for $\xi \in ]0,1[$.

By choosing $x_{k-1}, x_k$ as nodes, the secant method is obtained:
\begin{equation}\label{sec}
\Delta_k=-f[x_{k-1},x_k]^{-1}f(x_k)
\end{equation}
\hskip 12pt This is a superlinear method (its order of convergence is $(1+\sqrt{5})/2$), but the derivatives are not explicitly needed, and it only requires one additional evaluation of $f(x)$ for every iteration.

Finally, our goal in this paper is addressed to Steffensen's method. It is a quasi-Newton method in which the nodes are $x$ and $x+f(x)$:
\begin{equation}\label{stef}
\Delta_k=-f[x_k+f(x_k),x_k]^{-1}f(x_k)
\end{equation}
\hskip 12pt This method is second order, as Newton's, and it requires only two evaluations of $f(x)$ but the derivative is not needed. However, in spite of the amount of literature it has originated (as in \cite{jain,cand,chen,johnson,michel}), it is not very popular and it is often dismissed as an alternative to other classical methods. The reasons, even  though not explicit, are related to dimensionality.

As opposed to the other methods we have outlined in this section (except perhaps (\ref{firstorder})), Steffensen's is scale dependent. The sum $x+f(x)$ depends strongly on the scales of both $x$ and $f(x)$. A change of scale on $x$, $cx$, or $f(x)$, $kf(x)$, transform $x+f(x)$ in $cx+kf(x)$ which is not a change of scale in the node.

Some variations have been proposed in order to adjust this method. One of the most rediscovered is the modification of the node $x+f(x)$ by a constant $\alpha$ multiplying the function: $x+\alpha f(x)$. The idea behind this is the same useful resource after many scaling problems: a small value of $|\alpha|$ provides a better approximation of $f[x+\alpha f(x),x]$ to the exact derivative $f'(x)$ (see \cite{cand}).

While being a good alternative to the method, the problem is deeper. It is not a scaling  but a dimensional trouble: $[x]$ and $[f]$ are usually different and their sum is inconsistent.

Part of the problem can be solved by considering $\alpha$ not as an adimensional factor but as $[\alpha]=[x][f]^{-1}$ (that is, with the same dimension as $f'(x)^{-1}$). In this case, difficulties may arise from the own definitions of $\alpha$, which can vary largely if so does the derivative. A possibility is to use the same idea we did in (\ref{alter}):
$$
f[x,x+\lambda f(x)/f'(x_0)]
$$
\hskip 12pt Here, $\lambda$ is an adimensional constant and it is not influenced by a change of scale.

\vskip 6pt

Summing up, bisection, Newton's and secant methods are robust due to their scale independence, while the first order one (\ref{firstorder}) and Steffensen's are not and some adjustments must be done in order to keep them operative.

The emphasis on the importance of Newton's  is not a personal appreciation. The following theorem shows a dimensional aspect of this method which deserves to be highlighted:


\begin{theorem}{}\label{th1}
Under the notation in this paragraph, if $f(x)$ is once differentiable, then, any iterative method in the form (\ref{delta}) must verify the following dimensional relation:
$$
[\Delta_k]=\left[\frac{f(x_k)}{f'(x_k)}\right]
$$
As a consequence, any iterative method can be represented as:
\begin{equation}\label{carac}
\Delta_k=-h(x_k)\frac{f(x_k)}{f'(x_k)}
\end{equation}
\noindent where $h(x)$ is an adimensional function, $[h]=[1]$.
\end{theorem}


The proof is evident, because $[\Delta_k]=[x]=[f/f']$. However, some of its consequences are not: for instance, it is not so evident to think of bisection as an adimensional  modification of Newton's.

Moreover, by taking into account that $[f/f']=[x]$, its derivative 
$$
\left(\frac{f(x)}{f'(x)}\right)'=1-\frac{f''(x)f(x)}{f'(x)^2}=1-L_f(x)
$$
\noindent is adimensional. $L_f(x)$ is called the {\it degree of logarithmic convexity}, and (\ref{carac}) can be expressed as:

\begin{equation}\label{carac2}
\Delta_k=-h(L_f(x_k))\frac{f(x_k)}{f'(x_k)}
\end{equation}
\hskip 12pt By differentiating (\ref{carac2}) and cancelling derivatives, we obtain this well known characterization of second and third order algorithms \cite{gander}:

\begin{corollary}
An iterative method as (\ref{carac2}) is second order if $h(0)=1$, and third order if, besides, $h'(0)=1/2$.
\end{corollary}

From {\it Theorem 1}, higher order characterizations can be obtained, depending on higher derivatives of $L_f(x)$.

\section{Adimensional polynomials and functions}\label{sec3}

Section \S\ref{sec2} provides some insight on the relation between iterative methods, scaling and dimension. How do linear changes in $x$ or $f(x)$ affect the iterations?

As we stated above, Newton and secant are independent of the scale: 

If $y=cx$, and $g(y)=kf(y)$, and $y_k=cx_k$, then 
$$
y_{k+1}=y_k-\frac{g(y_k)}{g'(y_k)}=c\left(x_k-\frac{f(x_k)}{f'(x_k)}\right)=cx_{k+1}
$$

Rescaling $f(x)$ does not modify the iteration, in the sense that a change of scale in the variable $x$ produces the same rescaling on the iterates.

Similar behavior shows the secant method and even bisection. But first order (\ref{firstorder}), and Steffensen's are different. Scale problems in the first one can be avoided by a slight redefinition of the variable, denoting $\alpha=\lambda/K_1$
\noindent where $K_1=\sup\{f'(z): z\in U\}$, and $U$ the domain where the iteration is defined. In this case, a rescaling of $x$ or $f(x)$ modifies $K_1$, but not $\lambda$.

In order to analyze Steffensen's method, some preliminary comments must be done on the estimation of the errors of iterative methods.

There are two principles to ensure the convergence of the methods in Banach spaces. One of them is the majorizing sequence and the other one is  the majorizing function. 

\begin{definition} Given a sequence $\{x_n\}_{n\geq0} \subset E$ ($E$ a Banach space), an increasing sequence $\{t_n\}_{n\geq0} \subset {\mathbb R}$ is a majorizing sequence if, for all $k\geq0$: 
$$
||x_{k+1}-x_k|| \leq t_{k+1}-t_k
$$
\end{definition}

Convergence of $\{t_n\}_n$ implies convergence of $\{x_n\}_n$ and, if $t^*$ and $x^*$ are their limits, $||x_k-x^*||\leq t^*-t_k$.

Majorizing sequences are often obtained from application of the iterative method to majorizing functions, which are defined as:

\begin{definition} Given $F: E \longrightarrow V$, a real function $f(t)$ is a majorizing function of $F(x)$ if for an $x_0\in E$, for any $x$ in a neighborhood of $x_0$, $||F(x)|| \leq |f(||x-x_0||)|$.
\end{definition}

Newton's method ($x_{n+1}=x_n-F'(x_n)^{-1}F(x_n)$) illustrates this semilocal convergence analysis. Kantorovich established conditions on the initial guess, $x_0$, in order to get convergence of Newton's. These conditions were related to the boundness of the following variables:

\begin{theorem}\label{newton} Let $F:U \longrightarrow V$ be a differentiable function defined in a neighborhood $U\subset E$, and $x_0 \in U$ such that $F'(x_0)$  is invertible, and there exist three positive constants $K_2,\eta, B>0$ verifying:

\begin{enumerate}
\item $||F''(x)|| \leq K_2$, for all $x \in U$.
\item $||F'(x_0)^{-1}|| \leq B$.
\item $||F'(x_0)^{-1}F(x_0)|| \leq \eta$
\end{enumerate}

Then, if $K_2B \eta \leq 1/2$, Newton's method from $x_0$ converges.
\end{theorem}

The proof can be found in elementary textbooks \cite{ost,ort}. For our purposes, we point out a key step in this proof: the polynomial
$$
p(t)=\frac{K_2}{2} t^2-\frac{1}{B}t+\frac{\eta}{B}
$$
\noindent is a majorizing function, when starting from $t_0=0$. Condition $K_2B\eta \leq \frac{1}{2}$ means that the discriminant of $p(t)$ is greater or equal to zero and, therefore, $p(t)$ has two simple positive roots or a double positive root.

Dimensional analysis simplifies this polynomial. We consider $a=K_2B\eta$, $s=t/\eta$ and we define $q(s)$:
$$
p(t)=\frac{B}{\eta} q(s)
$$
\hskip 12pt From this definition, $q(s)=(a/2)s^2 -s +1$ (the minus sign in the linear term is required in order to get an increasing sequence $\{t_k\}_{k\geq 0}$). 

$q(s)$ is an adimensional version of $p(t)$ invariant from scales. If there is a change $\widetilde{t}=ct$ and $\widetilde{p}=kp$, then their correponding $\widetilde{K}_2,\widetilde{\eta}$, $\widetilde{B}$ are linearly modified, but the adimensional polynomial $q(s)$ remains ($\widetilde{a}=\widetilde{K}_2 \widetilde{B} \widetilde{\eta}=K_2 B \eta=a$).

We call $q(s)$ the {\bf adimensional form} of $p(t)$. From an scaling point of view, an adimensional form is a normalization of the polynomial such that $q(0)=1$ and $q'(0)=-1$.

The concept of adimensional form can be generalized to any nonlinear function such that both $f(0)$ and $f'(0)$ are not equal to zero. The adimensional form of a third degree polynomial:
$$
p(t)=\frac{K_3}{6}t^3 + \frac{K_2}{2}t^2-\frac{1}{B}t+\frac{\eta}{B}
$$
is: $q(s)=(b/6)s^3+(a/2)s^2-s+1$.

\noindent $b=K_3B\eta^2$, $a=K_2B\eta$. This third order polynomial is used for obtaining conditions on third order methods, such as Halley's, Chebyshev's or Euler's \cite{nos,nos2}. Convergence is obtained, as for Newton's,  when $q(s)$ has two simple positive roots or one positive double root (there is always a negative root).

Further generalization leads to the following:

\begin{definition} Given a nonlinear function $f(x)$, such that both $f(0)\neq 0$ and $f'(0)\neq 0$, the {\bf adimensional form} of $f(x)$ is a function $g(y)$ such that $g(0)=1$, $g'(0)=-1$ and $y=-(f'(0)/f(0))x$.
\end{definition}

Then, given $f(x)$ such that $f(0)\neq 0$ and $f'(0)\neq 0$, 
$$
f(x)=f(x_0)\left(g\left(-\frac{f'(0)}{f(0)}x\right)\right)
$$
\hskip 12pt We notice that $[f(x_0)][g]=[f]$, and therefore $[g]=[1]$, and $[y]=[f'][f]^{-1}[x]=[1]$. So, both $y$ and $g(y)$ are adimensional and they are scale independent, which is the goal of this definition.

\vskip 6pt

In \cite{nos,nos2}, a different way to study convergence is introduced. The idea is to get a system of sequences, two of them fundamental because they control both the iterates $x_k$ and their values $F(x_k)$, and some auxiliary sequences controlling different aspects of the iteration, if needed. The main advantage of these systems is that the sequences are adimensional, and the error estimates do not depend on the scales. We will use that technique in the next section applying it to Steffensen's method.

\section{A system of {\it a priori} error estimates for Newton's method}\label{sec4}

We start by illustrating the use of systems of bounds in the last paragraph, applying it to Newton's method. We refer the interested reader to \cite{nos} for details. Let us consider the following system of sequences:

Let $a>0$, $a_0=1$, $d_0=1$. For all $n \geq 0$. By recurrence, we define the following sequences
\begin{align}\label{sist}
{\rm {\bf (i)}} \ & a_{n+1}=\frac{a_n}{1-aa_nd_n}   &
{\rm {\bf (ii)}} \ & d_{n+1}=\frac{a}{2}a_{n+1}d_n^2 
\end{align}
\hskip 12pt This system is said to be {\bf positive} if $a_n>0$ for all $n \geq 0$.
It is easy to prove that the system is positive if and only if $0<a<1/2$. Positivity and convergence of $\sum_{n\geq 0} d_n$ are equivalent.

It can be checked that, for all $n\geq 0$:
\begin{equation}\label{cons}
\left(\frac{1}{a_n}\right)^2-a\frac{d_n}{a_n}=1-2a
\end{equation}
\hskip 12pt Therefore, $a_n$ depends on $d_n$, and:
\begin{equation}\label{dn}
d_{n+1}=\frac{a}{2}\frac{d_n^2}{ad_n+\sqrt{(ad_n)^2+(1-2a)}}
\end{equation}\label{dn}
(This function relating $d_n$ and $d_{n+1}$ is called {\it rate of convergence} in \cite{ptak}).

Thus: 
\begin{equation}\label{sum}
\sum_{n=0}^{+\infty} d_n=\frac{1-\sqrt{1-2a}}{a}=:s^* 
\end{equation}
Now, this theorem gives semilocal conditions for the convergence of Newton's method in Banach spaces ($(E,|| \cdot ||)$, $(V, || \cdot ||)$):

\begin{theorem}\label{q} Given a function $F: E \longrightarrow V$ at least twice differentiable, $x_0 \in E$, such that $F'(x_0)$ is invertible, and there exists positive constants $K_2,B,\eta$ verifying:

\begin{align*}
{\rm {\bf (i)}} \ & a=K_2B\eta \leq 1/2\\
{\rm {\bf (ii)}} \ & ||F''(x)|| \leq K_2 \quad \text{for all} \ x \in E \\
{\rm {\bf (iii)}} \ & ||F'(x_0)^{-1}|| \leq B \\
{\rm {\bf (iv)}} \ &  ||F'(x_0)^{-1}F(x_0)|| \leq \eta
\end{align*}
Then, Newton's method ($x_{n+1}=x_n-F'(x_n)^{-1}F(x_n)$) converges to a root, $x^*$, of $F(x)$ and the following inequalities hold:
\begin{enumerate}
\item $||F'(x_n)^{-1}|| \leq a_nB$, for all $n \geq 0$.

\item $||x_{n+1}-x_n|| \leq d_n \eta$, for all $n \geq 0$.

\item $||x^* -x_0|| \leq s^* \eta$.

\item $||x^*-x_n|| \leq (s^* -\sum_{k=0}^{n-1}d_k) \eta$, for all $n \geq 1$.

\end{enumerate}

\noindent $\{a_n,d_n\}$ is the system generated by (\ref{sist}) from $a$.
Furthermore, $x^*$ is the only root of $F(x)$ in the ball $\mathcal{B}(x_0,s^{**}\eta)$, where  
$$
s^{**}:=\frac{1+\sqrt{1-2a}}{a}
$$
\end{theorem}

The proof of this theorem is a direct consequence of the adimensional polynomial $q(d)=(a/2)s^2-s+1$ related to the majorizing one: $p(t)=(K_2/2)t^2-t/B+\eta/B$. Under the hypothesis of the theorem, $q(s)$ has two positive roots, $s^* \leq s^{**}$, and Newton's sequence: $s_0=0$, $s_{n+1}=s_n-q(s_n)/q'(s_n)$ converges monotonically to $s^*$ and verifies: $s_{n+1}-s_n=d_n$ and $-1/q'(s_n)=a_n$.

The adimensional polynomial explains also why (\ref{cons}) and (\ref{sum}) verify: the value $q'(s)^2-4q''(s)q(s)$ is an invariant for any quadratic polynomial, and its value is the discriminant, which in this case is $1-2a$. And $s^*=\lim \{s_n\}=\sum_n d_n$.

On the other side, it is easy to check that, if $a \neq 0$:
$$
\sum_{k=0}^{n-1}d_k=\frac{1}{a}(1-\frac{1}{a_n})
$$
\noindent given that $-1/a_n=q'(s_n)=as_n-1$, and $s_n=\sum_{k=0}^{n-1} d_k$.

We finish this section by introducing a new definition of order of convergence based on the adimensional functions and sequences.

\vskip 6pt

\begin{definition} Let us assume $(E,|| \cdot ||)$ is  Banach space, and $\{x_n\}_n \subset E$. Then, $p$ is the {\bf adimensional order} of $\{x_n\}_n$ if there exists an adimensional scalar sequence $\{d_n\}_n$ with order $p$, and a constant $\eta$ with $[\eta]=[x_0]$ such that:
$$
 ||x_{n+1}-x_n|| \leq d_n \eta
$$
 
\end{definition}
As expected, R-order and adimensional R-order  (AR-order) are equivalents (R-order itself is always adimensional), but adimensional Q-order (AQ-order) lies between classic Q and R-order (Q-order implies AQ-order, and AQ-order implies R-order). In fact, AQ-order behaves as classic Q-order but applied to relative instead of absolute errors.

From (\ref{dn}) and {\bf 2.}, it is deduced that Newton's method has second AR-order when the root is simple ($a<1/2$), and first AR-order when $a=1/2$. Actually, the discriminant $1-2a$ marks the boundary between fast and slow convergence \cite{ptak}.

\vskip 6pt

This is an example illustrating how adimensional polynomials and dimensional analysis make easier the study of the iterative method. In the next paragraph we go further in order to devise robust methods.

\section{Steffensen's method}\label{sec5}

As it was defined in the introduction, Steffensen's method for a scalar function $f(x)$ is defined, from $x_0$, by the iteration:
\begin{equation}\label{stefsc}
x_{n+1}=x_n-f[x_n+f(x_n),x_n]^{-1}f(x_n) \quad \text{for all} \ n \geq 0
\end{equation}
\hskip 12pt The divided difference 
$$
f[x,y]=f[y,x]=\frac{f(y)-f(x)}{y-x}
$$
\noindent is a well known operator in interpolation. We point out just one property arising when $f(x)$ is twice differentiable: there exists a value $\xi$  ($x< \xi<y$, or $y<\xi<x$) such that:
\begin{equation}\label{aprox}
f[x,y]=f'(x)+\frac{f''(\xi)}{2}(y-x)
\end{equation}
\hskip 12pt Thus, divided differences approximate the derivative of the function. A good choice of the nodes is a base to obtain an approximation to Newton's method when derivatives are unknown or costly. One of the simplest of these methods is the secant one, where the chosen nodes are the two previous iterates:
$$
x_{n+1}=x_n-f[x_{n-1},x_n]^{-1}f(x_n)
$$
\hskip 12pt In spite of a small reduction of the order of convergence, and an increase in the memory ($x_{n-1}$ must be stored), only one function evaluation is needed. This is a reason for being a very popular method. Let us notice that the dimension of these $f[x_{n-1},x_n]$ is $[f][x]^{-1}=[f']$. So, secant method does not introduce issues on scales or dimensions.

\vskip 6pt

After all these preliminaries, we deal with Steffensen's. Here, we consider the nodes $x_n$ and $x_n+f(x_n)$. The key idea is that, if the root $x^*$ is simple, by the mean value theorem, $f(x_n)=O(x_n-x^*)$, and therefore, $f[x_n,x_n+f(x_n)]$ is a better approximation to $f'(x_n)$ than that of the secant method. In fact, when it converges to a simple root, it is a second order method. Besides sharing the fact that no derivative is needed, Steffensen's has the advantage of a greater convergence order and no additional storage is required, unlike it happens to secant.

As we said in the introduction, a great amount of papers have appeared related to Steffensen's method, but it is not widely used because its scale dependence and dimension inconsistence. From a dimensional point of view, $[x]$ and $[f(x)]$ are different, and it is not a good idea to add both terms. Even in those cases where $[x]=[f]$, a different rescaling of the variable $cx$ and the function $kf$ makes the sum $cx+kf(x)$ unrelated to the original $x+f(x)$. 

Numerical analysts provide a sort of solution to this problem, based rather on scales than dimensions. Instead of evaluating $x+f(x)$, a parameter $\alpha$ is chosen in order to improve the method. The smaller $|\alpha|$ is, the better $f[x,x+\alpha f(x)]$ approximates $f'(x)$. Nevertheless, in order to keep consistence, $\alpha$ must verify $[\alpha]=[x][f]^{-1}$. A definition of the type: $\alpha=\lambda/f'(x_0)$ makes $
f[x,x+\lambda f(x)/f'(x_0)]$ independent from rescaling, and $\lambda$ is, in this case, adimensional, as searched. On the other side, as $[y]=[|y|]$, $f[x,x+\lambda f(x)/|f'(x_0)|]$ can be easier to implement, especially, as we will see, when $x$ and $f(x)$ are vectors instead of scalars. Of course, the choice of $f'(x_0)$ is arbitrary, and any other value of $f'$, or any variable with the same dimension as $f'$  can be used.

Adimensional functions provide a different point of view, which will be explained in the rest of this paper. We begin by finding a system of bounds for Steffensen's, similar to that of Newton's.

\subsection{A system of {\it a priori} error estimates for Steffensen's method}

Let's try to understand the behaviour of Steffensen's method with a simple example. We consider the adimensional quadratic polynomial: $q(s)=(a/2)s^2-s+1$ and we remind that this polynomial does not have scale or dimensional issues. For these polynomials Steffensen's converges at least as fast as Newton's, and both are second order.

\begin{theorem} If  $0\leq a \leq 1/2$, $s_0=t_0=0$ and we denote the sequence of Steffensen's approximations:
$$
s_{n+1}=s_n-q[s_n, s_n+q(s_n)]^{-1}q(s_n)
$$
\hskip 12pt Then $\{s_n\}$ converges monotonically to $s^*$. Furthermore, if we denote Newton's approximations:
$$
t_{n+1}=t_n-\frac{q(t_n)}{q'(t_n)}
$$
it is verified that, for all $n\geq 0$, $0 \leq t_n \leq s_n \leq s^*$.
\end{theorem}

{\bf Proof:}

In the first place, we observe that Steffensen's iteration is well defined. For all $0 < s < s^*$, 

$$
q[s,s+q(s)]=\frac{q(s+q(s))-q(s)}{q(s)}=q'(s)+\frac{a}{2}q(s)=:g(s)
$$

$g(s)$ is an increasing function in $[0,s^*]$, 

$$
g'(s)=q''(s)+\frac{a}{2}q'(s)=a+\frac{a}{2}(as-1)=\frac{a}{2}(1+as)>0
$$

And, $-1+a/2=g(0) < g(s) < g(s^*)=(a/2) q'(s^*) \leq 0$.

Therefore, $q[s,s+q(s)] \neq 0$ and it is invertible.

On the other side, denoting by $\Delta_S$ the Steffensen's correction and by $\Delta_N$ the Newton's one, we have:
\begin{equation}\label{compare}
\Delta_S(s)=-q[s,s+q(s)]^{-1}q(s)=-\frac{q(s)}{q'(s)+\frac{a}{2}q(s)} \geq -\frac{q(s)}{q'(s)}=\Delta_N(s)
\end{equation}
$$
\Delta_S(s_0)=\Delta_S(0)=\frac{1}{1-\frac{a}{2}}\leq \frac{1-\sqrt{1-2a}}{a}=s^*=s^*-s_0
$$
$$
\Delta_S(s^*)=0
$$
\hskip 12pt Differentiating $\Delta_S(s)$ with respect to $s$,
$$
\Delta_S'(s)=-\frac{q'(s)^2-aq(s)}{(q'(s)+\frac{a}{2}q(s))^2}
$$
\hskip 12pt As
$$
q'(s)^2-aq(s)=(q'(s)^2-2aq(s))+aq(s)=(1-2a)+aq(s)\geq 0
$$ 

$\Delta_S(s)$ is decreasing and positive, and, if $0\leq s \leq s^*$, then $s+\Delta_S(s)\leq s^*$.

As $\Delta_S(s) \geq 0$, the sequence $\{s_n\}_n$ is increasing and bounded, and, hence, it converges. The limit is, obviously, $s^*$.

On the other side, monotony and (\ref{compare}) prove that:
$0\leq t_n \leq s_n$ for all $n \geq 0$.
\hfill $\blacksquare$

\vskip 6pt

However, if we take $\eta > 2/a$, and $K_2=a/\eta$,
the polynomial $p(t)=(K_2/2)t^2-t+\eta$ has the same adimensional form $q(s)=(a/2)s^2-s+1$, but $p[0,p(0)]>0$ and, if $t_0=0$, then $t_1<t_0$ and the sequence of Steffensen's iterates does not converge. If $\eta=2/a$ it gets worse, because $p[0,p(0)]=0$ and, as it is not invertible, the sequence is not even defined.

These phenomena do not happen in scale and dimensional independent methods such as Newton's or secant, but they may become a real pain in Steffensen's.

We construct now a system of estimates for Steffensen's applied to adimensional polynomials. Given a positive $a\geq 0$ and $c_0=a_0=1$, $r_0=0$, we define for all $n \geq 0$:
\begin{align*}
{\bf (i)}  & \ b_n=\frac{a_n}{1-\frac{a}{2}a_nc_n} \\
{\bf (ii)} & \ d_n=b_nc_n \\
{\bf (iii)} & \ a_{n+1}=\frac{a_n}{1-aa_nd_n} \\
{\bf (iv)} & \ c_{n+1}=\frac{a^2}{2}d_n^2(r_n+\frac{1}{2}c_n)=\left\{ \left(\frac{1}{a_n}+1\right)+\frac{a}{2}c_n \right\} \frac{a}{2}d_n^2 \\
{\bf (v)} & \ r_{n+1}=r_n+d_n
\end{align*}
Now, our system consists of five sequences, two of them fundamental, $\{a_n\}_n$ and $\{c_n\}_n$, and the other three are auxiliary just for simplification.

The following theorem holds:

\begin{theorem}\label{q} Given $q(s)$, adimensional polynomial, and the sequence $\{s_n\}_{n\geq 0}$ generated by $s_0=0$, $s_{n+1}=s_n-q[s_n,s_n+q(s_n)]$, $n\geq 0$, then:

\begin{align*}
{\bf (I)} \ &  |q'(s_n)^{-1}|=a_n. \\
{\bf (II)} \ & |q[s_n, s_n+q(s_n)]^{-1}|=b_n. \\
{\bf (III)} \ & |q(s_n)|=c_n. \\
{\bf (IV)} \ & |s_{n+1}-s_n|=d_n. \\
{\bf (V)} \ & s_n=r_n. 
\end{align*}
\end{theorem}

{\bf Proof:} By induction. If {\bf (III)} and {\bf (IV)} are valid for $n$, then:
\begin{align*}
& {\text {\bf (I)}} & \ q'(s_{n+1})=q'(s_n)+a(s_{n+1}-s_n) \ {\rm implies} 
\end{align*}
$$
\frac{1}{q'(s_{n+1})}=\frac{1}{q'(s_n)+a(s_{n+1}-s_n)}=\frac{1}{\frac{-1}{a_n}+ad_n}=\frac{-a_n}{1-aa_nd_n}=-a_{n+1}
$$
\begin{align*}
& {\text {\bf (II)} }  & \  q[s_n, s_n+q(s_n)]=q'(s_n)+\frac{a}{2}q(s_n)=\frac{-1}{a_n}+\frac{a}{2}c_n
\end{align*}
Then,
$$
q[s_n, s_n+q(s_n)]^{-1}=\frac{-a_n}{1-\frac{a}{2}a_nc_n}=-b_n
$$
\begin{align*}
& {\text {\bf (III)}} & \  q(s_{n+1})=q(s_n)+q'(s_n)(s_{n+1}-s_n)+\frac{a}{2}(s_{n+1}-s_n)^2
\end{align*}
The following expression must be $0$:
\begin{equation}\label{c}
q(s_n)+q[s_n, s_n+q(s_n)](s_{n+1}-s_n)=q(s_n)+q'(s_n)(s_{n+1}-s_n)+\frac{a}{2}q(s_n)(s_{n+1}-s_n)
\end{equation}
We have:
\begin{equation}\label{c2}
q(s_{n+1})=\frac{a}{2}d_n(d_n-c_n)
\end{equation}

From (\ref{c}), $q(s_n)=-q'(s_n)d_n-(a/2)q(s_n)d_n$

And  (\ref{c2}) becomes:
$$
q(s_{n+1})=\frac{a}{2}d_n^2 (1+q'(s_n)+\frac{a}{2}q(s_n))=
-\frac{a}{2}d_n^2((q'(s_n)-q'(0))+\frac{a}{2}c_n)
$$
(because $q'(0)=-1$). Taking into account that $s_n=r_n$, finally:
$$
q(s_{n+1})=-\frac{a^2}{2}d_n^2(s_n+\frac{1}{2}c_n)=c_{n+1}
$$

The rest of the theorem is obvious.
\hfill $\blacksquare$

The system {\bf (I)}-{\bf (V)} shows that the sequence $\{d_n\}_n$ (similarly, $\{c_n\}_n$) is second AQ-order when $1-2a>0$, and first AQ-order when $2a=1$, as it happened wtih Newton. The rate of convergence, however, is not so simple. Some long and tedious but not so complicated algebraic manipulations, and the invariance of $(1/a_n^2)-2ac_n=1-2a$ gives:
$$
d_{n+1}=\frac{a}{2}\cdot \frac{(1+\alpha_n)+\frac{a}{2}\frac{\alpha_nd_n}{1+\frac{a}{2}d_n}}{\alpha_n\left(1-\frac{a}{2}\frac{d_n}{1+\frac{a}{2}d_n}\right)}\cdot d_n^2
$$
\begin{align*}
&{\rm with} & \ \alpha_n=\frac{ad_n}{1+\frac{a}{2}d_n}+\sqrt{\left(\frac{ad_n}{1+\frac{a}{2}d_n}\right)^2+(1-2a)}
\end{align*}

As in Newton, $1-2a$ is the boundary separating low and high speed of convergence. However, while this region of high order convergence is reached by $(ad_n)^2$ in Newton, in Steffensen it is reached by $(ad_n/(1+(a/2)d_n))^2$. This shows a faster (though with the same order) convergence of Steffensen's versus Newton's method.

\vskip 6pt

Now, we are ready to introduce the system of error bounds for Steffensen's method in Banach spaces. First, we must remind which kind of operator is the divided differences in Banach spaces. As opposed to scalar, there is not a unique divided differences operator for vectors. The characterization is based on the interpolation:

\begin{definition} Given $F: E \longrightarrow V$, $x,y \in E$, a linear operator $H:E \longrightarrow V$ is a divided difference of $F$ with nodes $x,y$ if $H(x-y)=F(x)-F(y)$. We denote $F[x,y]:=H$.
\end{definition}

If $F$ is twice differentiable, with $||F''(\xi)|| \leq K_2$ for every $\xi \in E$, it is verified:
$$
||F[x,y]-F'(x)|| \leq K_2 ||x-y||
$$
\hskip 12pt Further insight on divided differences in Banach spaces can be found in \cite{ptak}.  There, we can find one example of differentiable divided difference as follows:
$$
F[x,y]=\int_0^1 F'(x+\theta (y-x)){\rm d}\theta
$$
\hskip 12pt However, in practice, this integral is not explicitly evaluated, and divided differences are obtained from scalar divided differences.
 
The key point of the rest of this paragraph is to obtain an adimensional function $G$ related to $F$, which is the following.

From $x_0 \in E$, if there exists $F'(x_0)^{-1}$ and $F(x_0) \neq 0$, we get the change: 
\begin{equation}\label{y}
y=-\frac{F'(x_0)}{||F(x_0)||}x
\end{equation}
\hskip 12pt We define the adimensional form of functions in Banach spaces:

\begin{definition} Given $F$ under the assumptions of theorem (\ref{newton}), the {\bf adimensional form} of $F(x)$ is the function:
 \begin{equation}\label{adim}
G(y)=\frac{F(x)}{||F(x_0)||}
\end{equation}
where $y$ is defined in  (\ref{y}) and 
\begin{equation}\label{y0}
y_0=-\frac{F'(x_0)}{||F(x_0)||}x_0
\end{equation}
\hskip 12pt It is verified: $||G(y_0)||=1$ and $G'(y_0)=-I$.
\end{definition}

{\bf Remark:} Adimensional functions are not only useful for high order iterations. For example, the well known pathological zigzag optimization problem (\cite{strang}, p.348), consisting of minimizing the function $H(x,y)=(x^2+by^2)/2$, implying to solve, {\it via} steepest descent (a first order method similar to (\ref{firstorder})  the system: $F(x,y)=(F_1(x,y),F_2(x,y))=(x,by)=(0,0)$ presents a very slow convergence when $b$ is small. Actually, its rate of convergence is:
$$
H(x_{n+1},y_{n+1})=\left(\frac{1-b}{1+b}\right)^2 H(x_n,y_n)
$$
However, its adimensional form is $G(u,v)=-(x,y)/||F(x_0,y_0)||$, being $u=x/||F(x_0,y_0)||$, $v=b\cdot y/||F(x_0,y_0)||$, converges in just one step.

\vskip 6pt

Let us notice that the change of variables from $x$ to $y$ is also a change of space from $E$ to an adimensional space $\widetilde{V}$ homeomorphic to $V$. $\widetilde{V}$ is an example showing that dimension and dimensionality are different: it has the same dimension as $V$, but, for any $y\in \widetilde{V}$, $[y]=1$. As $G: \widetilde{V} \longrightarrow \widetilde{V}$ is an endomorphism, $y+G(y)$ is a consistent sum in $\widetilde{V}$, and it is possible evaluate $G[y,y+G(y)]$, which is independent from scales and dimensions, thus avoiding two of the main concerns about Steffensen's. We can consider $G$ as a preprocessing of the nonlinear operator and, certainly, its performance is similar as that of preprocessing of  matrices to solve linear systems: it reduces the computational cost after a small additional work before starting the method.

Instead of getting a root of $F$, the key point is to find a root of $G$. Any method is admissible but if the method is scale and dimensional independent, as Newton's, it is not worthy to transform the function. However, in Steffensen's, it may be very useful. For example, convergence is easier to obtain. Besides, the error bounds are more accurate, in the sense that they are optimal (optimality does not always guarantee better estimates, but  they ensure that under certain conditions they cannot be improved). Optimal bounds have not been found in general for this method. We will prove that the system of estimates as defined above provides optimal error bounds for Steffensen's method.

\subsection{Adimensional scale invariant Steffensen's method (ASIS method) and its convergence}

We call ASIS ({\it adimensional scale invariant Steffensen's}) to the Steffensen's method applied to the adimensional form of the function.

By the linear transformation from $F$ to $G$, as it was explained in (\ref{adim}), and $y_0$ in (\ref{y0}), we define for $n\geq 0$:
$$
y_{n+1}=y_n-G[y_n,y_n+G(y_n)]^{-1}G(y_n)
$$ 
$$
x_n=||F(x_0)||F'(x_0)^{-1}y_n
$$
\begin{theorem}\label{G}
Let us assume that there exists $G''$ in $V$, and $||G''(\xi)|| \leq a$, for $\xi \in V$ and $a\leq 1/2$. Then, from the system of bounds {\bf (i)-(v)} in the last subsection:
\begin{align*}
{\bf (a)} \ &  ||G'(y_n)^{-1}|| \leq a_n. \\
{\bf (b)} \ & ||G[y_n, y_n+G(y_n)]^{-1}|| \leq b_n. \\
{\bf (c)} \ & ||G(y_n)|| \leq c_n. \\
{\bf (d)} \ & ||y_{n+1}-y_n|| \leq d_n. \\
{\bf (e)} \ & ||y_{n+1}-y_0|| \leq r_n.
\end{align*}
\end{theorem}

{\bf Proof:} It is similar to that of theorem \ref{q}, by using Banach invertibility criterion. 

First, we note that $||G''(y)||\leq a$ for all $y$.

By induction, if {\bf (a)} and {\bf (c)} hold for $n=0$, then, for $n\geq 0$:
$$
||G'(y_n)-G[y_n,y_n+G(y_n)]||\leq \frac{a}{2} ||G(y_n)|| < \frac{1}{||G'(y_n)^{-1}||}
$$
(The last inequality by induction hypothesis). Then, by Banach invertibility criterion, $G[y_n,y_n+G(y_n)]$ is invertible and 
$$
||G[y_n,y_n+G(y_n)]^{-1}|| \leq \frac{||G'(y_n)^{-1}||}{1-\frac{a}{2}||G(y_n)||\cdot||G'(y_n)^{-1}||}
$$
\noindent and {\bf (b)} is verified. 

{\bf (d)}  and {\bf (e)} are evident.

By the same Banach criterion,
$$
||G'(y_{n+1})-G'(y_n)||\leq a ||y_{n+1}-y_n|| <\frac{1}{||G'(y_n)^{-1}||}
$$
\hskip 12pt Thus, $G'(y_{n+1})$ is invertible and {\bf (a)} is verified because:
$$
||G'(y_{n+1})^{-1}||\leq \frac{||G'(y_n)^{-1}||}{1-||G'(y_n)^{-1}||\cdot||y_{n+1}-y_n||}
$$
\hskip 12pt Finally, from:
$$
G(y_n)+G[y_n, y_n+G(y_n)](y_{n+1}-y_n)=0
$$
$$
G[y_n,y_n+G(y_n)]=G'(y_n)+G''(\xi)G(y_n) \ , \ \text{for}  \ \xi \in \widetilde{V}
$$
and
$$
||G(y_{n+1})-G(y_n)-G'(y_n)(y_{n+1}-y_n)|| \leq \frac{a}{2}||y_{n+1}-y_n||^2
$$
we have:
$$
||G(y_n)||\leq \frac{a}{2}||y_{n+1}-y_n||||(y_{n+1}-y_n)-G(y_n)(y_{n+1}-y_n)||
$$
$$
\leq \frac{a}{2}||y_{n+1}-y_n||^2||I-G(y_n)||\leq \frac{a}{2}d_n^2||G(y_n)+G'(y_0)||
$$
and it follows from 
$$
||G(y_n)+G'(y_0)||=||G'(y_0)-G'(y_n)-\frac{G''(\xi)}{2}(y_{n+1}-y_n)|
$$
$$
\leq \frac{a}{2}(a||y_n-y_0||+\frac{a}{2}||y_{n+1}-y_n||)
$$
\hfill $\blacksquare$

The estimates of this theorem are optimal, because they are reached for the adimensional polynomial $q(s)=(a/2)s^2-s+1$.

As a consequence of  theorem \ref{G}, the following main result proves the convergence of Steffensen's method in Banach spaces:

\begin{theorem}\label{prin} Let $F:E \longrightarrow V$ be a twice differentiable nonlinear function, $x_0 \in E$ such that $F(x_0) \neq 0$ and $F'(x_0)$ invertible. Let $K_2,B,\eta>0$ verifying $||F(x_0)||\leq \eta/B$, $||F'(x_0)^{-1}||\leq B$, and, for all $x\in E$, $||F''(x)|| \leq K_2$.

If $a=K_2B\eta \leq 1/2$, then the sequence generated by ASIS method converges to a root $y^*$, $G(y^*)=0$. Furthermore, if for any $n\geq 0$, $x_n=(\eta/B)F'(x_0)^{-1}y_n$, the following inequalities hold:
$$
||x_{n+1}-x_n|| \leq d_n \eta
$$
$$
||x^* -x_n|| \leq (s^* -s_n)\eta
$$
$$
F'(x_n) \ \text{is invertible and} \ ||F'(x_n)^{-1}|| \leq a_n B
$$
\end{theorem} 

\vskip 12pt

This theorem improves previous results. So far, sufficient conditions in literature implied the existence of a strictly positive constant $h>0$, depending on $K_2,B, \eta$ such that $a(1+h)<1/2$. These conditions are more restrictive than those of Newton and, hence, than those we propose in theorem \ref{prin}, where $a\leq 1/2$ suffices to ensure convergence.

We remark that not only invertibility but existence of $F[x_n,x_n+F(x_n)]$ is not ensured under these conditions; even $x_n+F(x_n)$ may not be evaluated. However, differences of  the adimensional form do exist and they are invertible.

\section{Numerical Examples}\label{sec6}

Our first example is a simple scalar equation:
$$
f_1(x)=\exp(x-1)-1=0
$$
with the obvious solution $x^*=1$.

We want to compare the performance of Newton's, Steffensen's and ASIS methods. Both $x$ and $f(x)$ are real. 

This is a case where Kantorovich conditions are more restrictive than the actual Newton iteration. Kantorovich ($a\leq 1/2$) guarantees convergence for $x_0 \geq 1-\log(\frac{1+\sqrt{3}}{2}) \approx 0.6881$. However, convexity of $f_1(x)$ enlarges the region of convergence, as we will see.

We start Newton's and Steffensen's from $x_0=0$. This initial guess does not verify Kantorovich conditions, but both methods converge. In the figure, we represent the logarithms of the errors, where it can be seen that Newton's converges faster than Steffensen's. Though both are second order, Steffensen's first iterations are clearly slower, because it takes some steps to {\it get into} the second order region. Steffensen's errors are not only greater than those of Newton's, but it requires some more iterations in order to get the computational accuracy.

The setting varies for ASIS. The adimensional version of $f(x)$ is 
\begin{equation}\label{adim1}
g(s)=\frac{1}{\exp(-1)-1}(\exp((e-1)s-1)-1)
\end{equation}
and $x=(e-1)s$.  Therefore, $s_0=0$. In figure 1 we represent the logarithms of the errors $|1-(e-1)s_n|$.
It can be seen that these errors are smaller than Newton's. That is, in this case, ASIS does not only improves Steffensen's, but Newton's, as  expected from the theoretical results above.

\begin{figure}[h!]
\centering
\begin{tabular}{cc}
\includegraphics[width=0.4\textwidth]{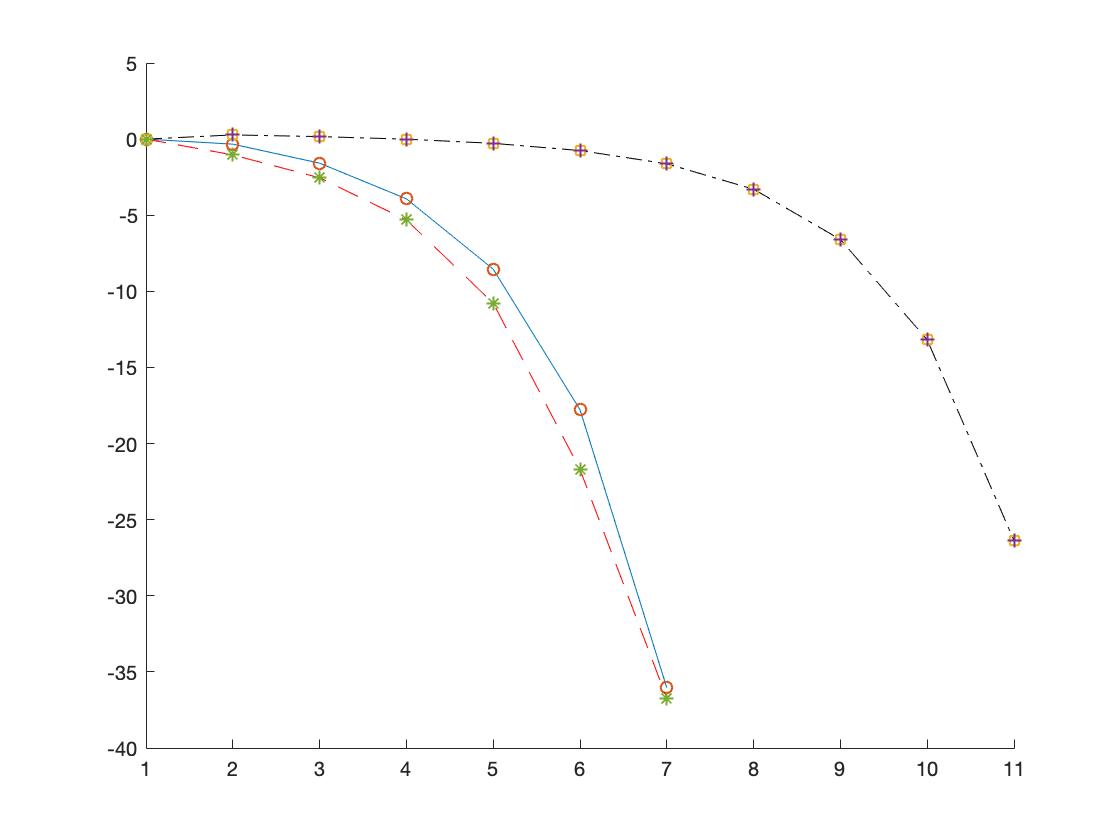} &
\includegraphics[width=0.4\textwidth]{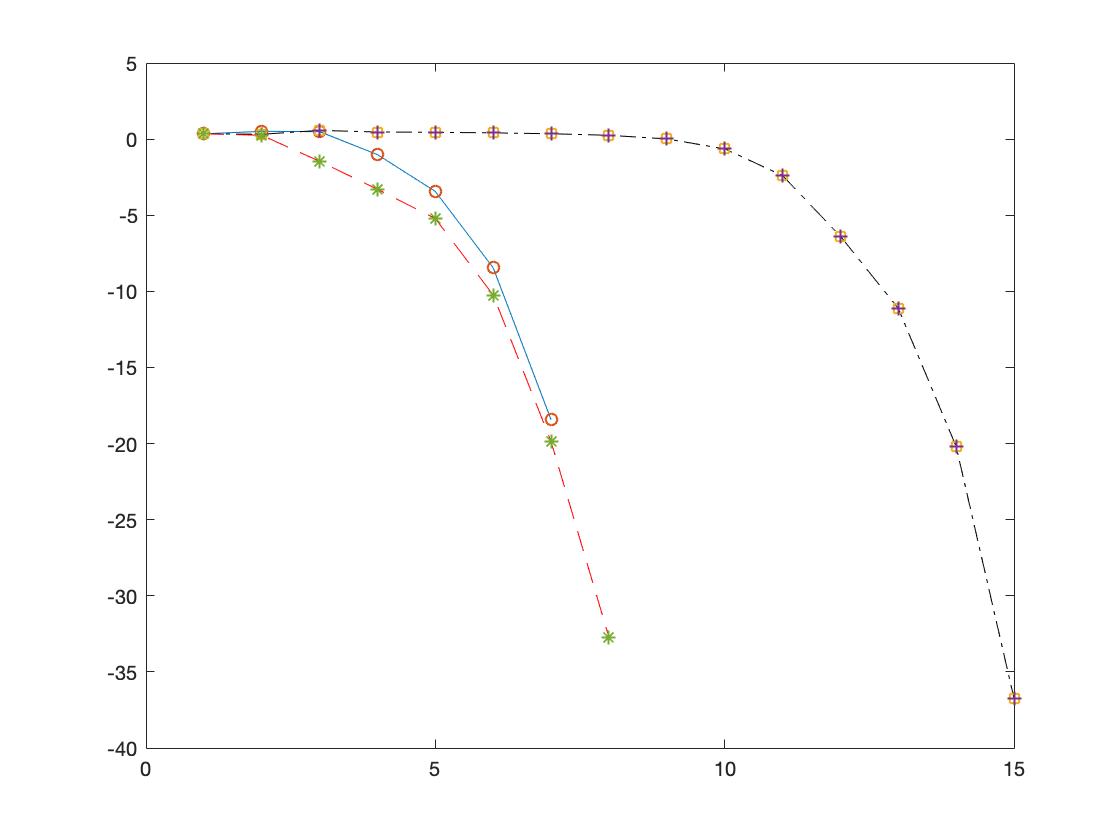}
\end{tabular}
\caption{Logarithms of errors: Newton, solid line, 'o'; Steffensen '.-.', '+'; ASIS '--', '*'.{\bf Example 1 (scalar)} (left); {\bf Example 3 (system)} (right)}
\end{figure}

A second example shows the consequences of rescaling the equation. We consider:
$$
f_2(x)=e^{2x-1}-1=0
$$
\hskip 12pt This function is obtained by: $f_2(x)=f_1(2x)$. The root of $f_2(x)$, $x^*=\frac{1}{2}$ is $\frac{1}{2}$ times the root of $f_1(x)$. Hence, Newton's errors for $f_2$ are one half the errors of $f_1$ due to the invariance respect to scales in the variable. The adimensional version of $f_2(x)$ is exactly the same as for $f_1(x)$, $g(s)$, and the errors of ASIS are also one half of the errors obtained for $f_1(x)$. As a consequence, ASIS and Newton's behave in the same way as in the first example (ASIS converges faster than Newton's).

However, Steffensen's applied to $f_2(x)$ needs $3705$ iterations to get an error less than $0.5$. We remark that $0.5$ is the initial error. Though it goes faster after that (iteration $3716$ provides an approximation less than $10^{-16}$, double precision accuracy), it is clearly an unpractical method in this case. In figure 2 we show the logarithms of the errors of classic Steffensen's (not adimensional): the first $3700$ iterations the errors look like a constant, and after that it behaves as a second order method. Its region of quadratic convergence is not so small, but it takes a very long way to reach it from the initial guess. It can get worse for finer rescalings: the method can blow up due to roundoff errors for $f(x)=f_1(cx)$ for large $c$. This behavior is usual for the usual Steffensen's method: its convergence and its speed of convergence, when it converges, depend strongly on the scaling of both, the variable and the own function.
\begin{figure}[h!]
\centering
\begin{tabular}{cc}
\includegraphics[width=0.4\textwidth]{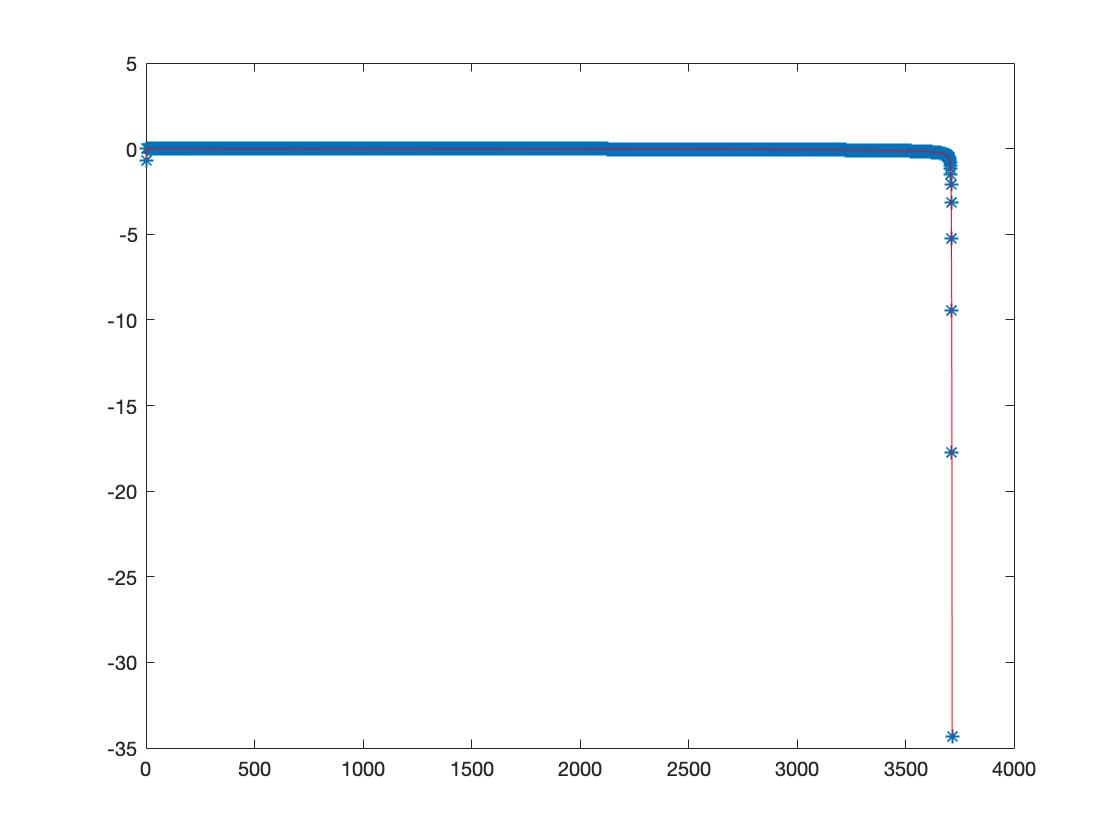} &
\includegraphics[width=0.4\textwidth]{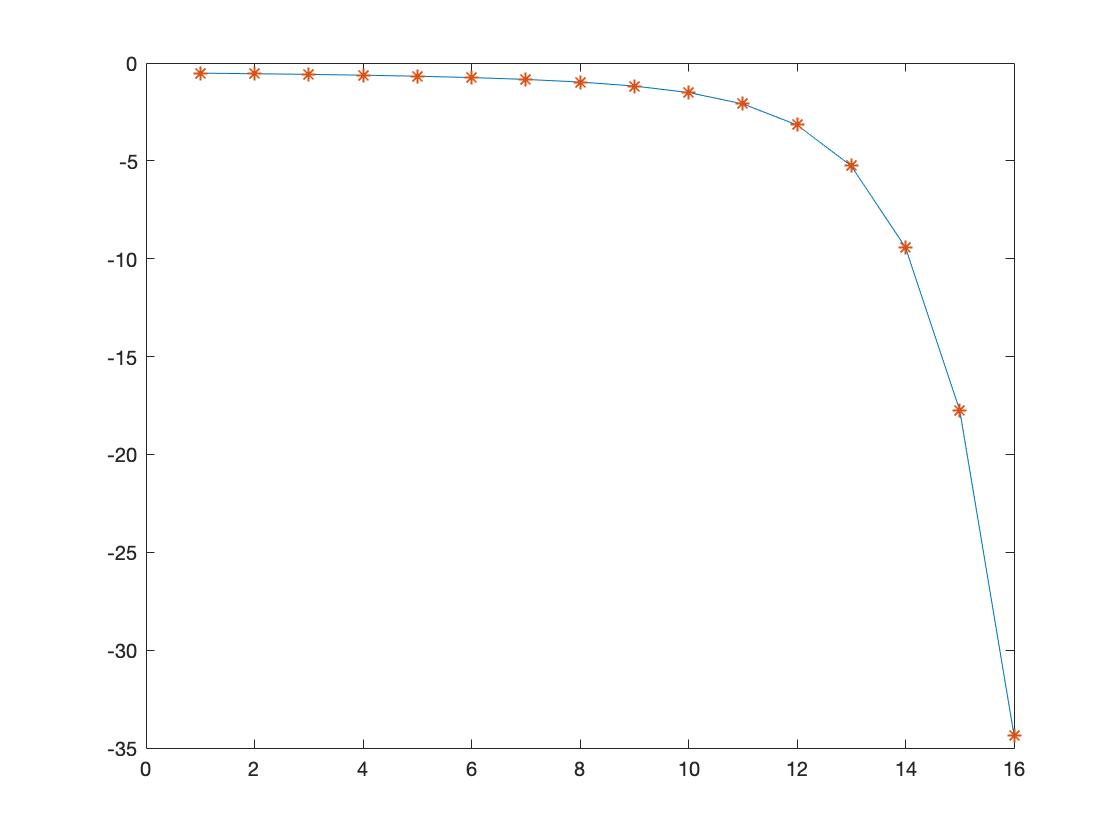}
\end{tabular}
\caption{Example 2: Logarithms of Steffensen's method for $f_2(x)$ (left); iterations from $3700$ to $3716$ (right).}
\end{figure} 
The last example in this section is a system of two nonlinear equations with two variables:
\begin{align*}
F_1(x,y)=&-4x(y-x^2+2)-2(1-x)=0\\
F_2(x,y)=&2(y-x^2+2)=0
\end{align*}
\hskip 12pt This system is obtained from an optimization problem, by unrestricted minimization of the function:
$$
F(x,y)=(y-x^2+2)^2+(1-x^2)^2
$$
\noindent which is the Euclidean norm of a quadratic curve.

Starting from $x_0=y_0=0$, both Newton and Steffensen converge. As it happens with our first one dimensional equation, Steffensen's is clearly slower than Newton. Among the possible choices of the divided difference operator, we used the one in \cite{ptak}.

Once again, we try the adimensional version of the system, normalized {\it via} the Euclidean norm, and we get:
\begin{align*}
g_1(s,t)=&-\frac{2}{3}s(\sqrt{5}t-\frac{5}{9}s^2+2)+(1+\frac{\sqrt{5}}{3}s)=0
\\
g_2(s,t)=&-(t-\frac{5}{9}s^2+\frac{2}{\sqrt{5}})=0
\end{align*}
\noindent with the change of variables:
$$
s=-\frac{3}{\sqrt{5}}x  \ ; \  t=\frac{1}{\sqrt{5}}y
$$
The fact that $\Delta F(0,0)$ is diagonal simplifies the change of variables (which becomes a diagonal transformation). On the other side, some good features of the original system, such as symmetry of the gradient, are lost in the adimensional system (therefore, the divided differences operator is no longer {\it quasi-symmetric}. 

Its behaviour is similar to that of the scalar case: ASIS is faster than Newton's, and  Steffensen's is clearly slower, as it can be seen in the figure 1.

\section{Conclusions}\label{sec7}

This paper was motivated by the study of quasi-Newton methods in optimization. It is somehow surprising that, in spite that Steffensen is a well known method (not as popular as Newton or secant, but well known) which can be implemented without explicit knowledge of derivatives, it is not used, and even proposed, as a feasible method in this topic. Some insight makes it clear that all the weakness of Steffensen's in scalar equations increase when dealing with equations in several variables. Nevertheless, most of these weaknesses are related to dimensional troubles. Independence from dimensions makes Steffensen's (and other methods) more robust and liable.

After some considerations about dimensionality, in this paper we introduced the adimensional form of functions and polynomials in order to analyze convergence of iterative methods. These adimensional forms  make it possible the design of methods based on classical ones that, adapted from scale dependent methods, overcome all concerns due to dimensionality and scaling. Here, we got the ASIS method, based on the classic Steffensen's method.

Semilocal conditions are established in order to ensure convergence of ASIS. ASIS method has the same convergence conditions as Newton's method. This is a really interesting result, because in general, conditions so far are more restrictive for classic Steffensen's than for Newton's. The estimates we settled are optimal in the sense that adimensional polynomials reach the estimates. This optimality is also an important improvement with respect to other estimators in literature, which are not optimal. Optimality improves the theoretical region of convergence. In fact, ASIS is at least as fast as Newton's for polynomials.

Semilocal conditions are not only stablished for theoretical knowledge, but for practical reasons: hybrid methods, using different algorithms, slower for the initial iterations, with the goal of getting approximations into the regions of convergence, and faster ones when the approximations are close enough to the exact solution, are often devised. In order to obtain good performance, well fitted estimators are needed. 

The price to pay in order to obtain adimensional functions is an increase of computational cost due to a change of variable which must be implemented. However, this cost can be reduced because in general $F'(x_0)^{-1}$ is not required to be obtained exactly, though the analysis of softening conditions is outside the scope of this paper. It is remarkable that even random numerical errors do have dimension. Randomness is a numerical condition, but it affects the measure of the parameters.

Finally, we remark that some conditions in the theorems in this paper can be relaxed. Not only the exact evaluation of the inverse of $F'(x_0)$, but also conditions on the second derivative, because it is enough to ask for Lipschitz continuity of the first derivative, or some other weaker conditions appearing often in literature. 

Last, but not least, ASIS is a good alternative method for quasi-Newton in optimization, because of its interpolatory conditions (no need to evaluate derivatives), robustness and speed of convergence.

The techniques developed here can be easily extended and generalized to other kind of algorithms. 

Summing up: though it is often disregraded, dimensional analysis is a powerful tool that helps to improve the performance of algorithms, and to a better understanding of their theoretical properties.  Adimensional functions get rid of drawbacks induced by heterogeneous data.

\vskip 6pt
{\bf Acknowledgements:}

First and mostly, this paper, and its author, are greatly indebted to the comments, suggestions and encouragement provided by Antonio Marquina, whose help improved both of them (paper and author).

\bibliographystyle{spmpsci} 

\end{document}